\documentclass[runningheads]{llncs}
\usepackage{graphicx}
\usepackage{amsmath, amsfonts, amssymb, makeidx}
\usepackage[cp1251]{inputenc}

\begin{document}
\title{Computational aspects of finding a solution asymptotics for a singularly perturbed system of differential equations\thanks{This research was performed in the framework of the state task in the field of scientific activity of the Ministry of Science and Higher Education of the Russian Federation, project "Development of the methodology and a software platform for the construction of digital twins, intellectual analysis and forecast of complex economic systems", grant no. FSSW-2020-0008.}}

\titlerunning{Computational aspects of finding solution asymptotics for system of PDEs}
%
\author{V.A. Krasikov\inst{1}\orcidID{0000-0001-9686-8493} \and A.V. Nesterov\inst{1}\orcidID{0000-0002-4702-4777}} 
\authorrunning{V.A. Krasikov, A.V. Nesterov}
%
\institute{Plekhanov Russian University of Economics \\
 Stremyanny Lane, 36, Moscow, 117997, Russia \\
\email{Krasikov.VA@rea.ru}\\
\email{Nesterov.AV@rea.ru}
}
\maketitle              
\begin{abstract}
We analyze the spatial structure of asymptotics of a solution to a singularly perturbed system of mass transfer equations. The leading term of the asymptotics is described by a parabolic equation with possibly degenerate spatial part. We prove a theorem that establishes a relationship between the degree of degeneracy and the numbers of equations in the system and spatial variables in some particular cases. The work hardly depends on the calculation of the eigenvalues of matrices that determine the spatial structure of the asymptotics by the means of computer algebra system Wolfram Mathematica. We put forward a hypothesis on the existence of the found connection for an arbitrary number of equations and spatial variables. 

\end{abstract}
%
%
%


\section{Introduction}

Using of the computer algebra systems for solving systems of differential equations or investigating properties of solutions by the means of symbolic computations is a very popular and fruitful approach~\cite{KrasikovJSFU2020,KrasikovLNCS2019}. One of the problems solved this way is symbolic research on eigenvalues and eigenvectors of matrices, associated with solutions of systems~\cite{DivakovTiutiunnikProgCompSoft,HsueIEEETransactions}. In the present work we investigate eigenvalues of matrices defining asymptotics of a solution to a specific system of differential equations by the means of Wolfram Mathematica.
		
In the article~\cite{b1} small parameter asymptotic expansion was built for a solution of the Cauchy problem for a singularly perturbed hyperbolic system of differential equations with multiple spatial variables
	\begin{equation}
   \label{c3}
      \varepsilon^{2} (U_{t} + \sum_{i=1}^{K}D_{i} U_{x_{i}}) = AU , \bar x\subset\Omega,t>0, 
\end{equation}
\begin{equation}
   \label{c4}
      U(\bar x, 0) = U^{0} (\bar x, \varepsilon), 
\end{equation}
 where $U(\bar x,t) = \lbrace u_{i} (\bar x, t)\rbrace$, \ $(i = \overline{1,n})  $  is a~vector of solutions, $0<\varepsilon <<1$ is a~small positive parameter, with a splash type initial conditions, concentrated in a small  $\varepsilon$-neighbourhood of the origin: 
	$$ 
	U(\bar x, 0) =  \omega (\bar x/\varepsilon)H,
	$$
		where $ H $ is a~vector, $\omega$ is 
		a~smooth function, rapidly decreasing as its argument goes to infinity along with all its derivatives.		
		   	
	Matrices~$D_i$ and~$A$ satisfy conditions described in~\cite{b1}, in particular, matrix $A$ has eigenvalue~$0$ with multiplicity~$1$  while real parts of other eigenvalues are negative. In the following 		
				 $h_{1}$ is the eigenvector of the matrix $A$ corresponding to the eigenvalue $\lambda_{1} = 0$; 
	$h_{1}^{*}$  is the eigenvector of the matrix $A^T$ corresponding to the eigenvalue $\lambda_{1} = 0$.

  If the conditions from the article~\cite{b1} are satisfied, asymptotic expansion of a~solution to the problem~(\ref{c3}) -~(\ref{c4}) for $t\geq t_0 >0 $ $\forall\; t_0$ independent of $\varepsilon $ may be represented in the form
\begin{equation}
\label{c6a}
	U(\bar x, t, \varepsilon)=\varphi_0(\bar \zeta,t)h_1+O(\varepsilon^2),
\end{equation}	
where $\bar\zeta=\{\zeta_i,i=1,...,N\}$,
	 
\begin{equation}
   \label{c6b}	
\zeta_i = \dfrac{x_{i} - v_{i}t}{\varepsilon},  v_{i} = (D_{i}h_{1}, h_{1}^*),
\end{equation}
and $\varphi_0(\bar \zeta,t)$ is a solution to the equation
\begin{equation}
   \label{c11}
       \varphi_{0t} + \sum_{i,j=1}^{K}M_{ij}\varphi_{0\zeta_{i}\zeta_{j}} = 0.
\end{equation}

Here \\
\begin{align}
	\label{m1}
	&M_{ii} = (\Psi_{i} G \Psi_{i} h_{1}, h_{1}^{*});\\
	\label{m2}
		&M_{ij} = ((\Psi_{i} G \Psi_{j} + \Psi_{j} G \Psi_{i}) h_{1}, h_{1}^{*})/2, i\ne j,
	\end{align}\\
where $\Psi_{i}=D_i-v_iE$, $E$ is the identity matrix, $G$ is a pseudo-inverse matrix for the matrix~$A$ (if~$(y,h^*_1)=0,$ a solution to the system of equations $Ax=y$  is given by $x=Gy+Ch_1$,  where $C$ is an arbitrary constant).
	The quadratic form	$\vartheta(\bar\zeta)=\sum_{i,j=1}^K M_{ij}\zeta_i\zeta_j$ satisfies the property of having fixed negative sign, which means that the equation (\ref{c11}) is parabolic.
	
	\section{Problem statement}
	
The main object of interest in the present work is a structure of the spatial part (which may be degenerate) of the parabolic equation~(\ref{c11}).
	To research a~possible degeneracy we have to find eigenvalues and eigenvectors of the matrix $M=\{M_{ij},1\leq i,j\leq K\}$.
	In the article~\cite{b2} for the system (\ref{c3})  with $n=2$ equations authors discovered the connection between the rank of the matrix~$M$ (the quadratic form~$\vartheta(\bar\zeta)$)  $\textrm{rank }M$   and the number of equations~$n$ in the system~(\ref{c3}): 
	$\textrm{rank }M=1$ if $n=2$ , $\textrm{rank }M=2$  if $n\geq 3.$ 
			
	The system~(\ref{c3}) with greater number of equations was considered in the article~\cite{b1}. Numeric computations of eigenvalues and eigenvectors (for given number matrices~$D_i,A$)	lead to the hypothesis on connection of the quadratic form~$\vartheta(\bar\zeta)$ rank  and the number of spatial variables~$K$ and the number of equations~$n$ in the system~(\ref{c3}).
	

In the present work we obtained the confirmation on the proposed in~\cite{b1} hypothesis for greater number of spatial variables~$K$ and greater number of equations~$n$ of the system by the means of symbolic computations of eigenvalues and eigenvectors of the matrix~$M$ using the computer algebra systems. 
The main result is

\begin{theorem}\rm

 Let the matrix~$M$ be defined by the formulas~ (\ref{m1}). Then its rank
  is equal to $$\textrm{rank }M=\left\{\begin{array}{l}n-1, \textrm{ if }n\leq K, \\ K, \textrm{ otherwise,} \end{array}\right.$$
where $n=2,\ldots,5, K=2,\ldots,5$ are the number of equations and the number of spatial variables in the system~(\ref{c3}) respectively.
\label{th:main}
\end{theorem}

Theorem~\ref{th:main} is proved for $n=2,\ldots,5, K=2,\ldots,5,$ the proof for arbitrary $n, K$ is under development.
\vskip0.4cm
\noindent{\bf Hypothesis 2.}
The statement of Theorem~\ref{th:main} holds for any $n\geq 2, K\geq 2.$

 \section{The computation of the eigenvalues of the matrix~$M$}   

		To prove Theorem~\ref{th:main} we found the analytic expressions for the eigenvalues and the eigenvectors of~$M,$
			which was computed by using the initial data -- namely, the $D_i,A$ matrices by the means of the formulas~(\ref{m1}).
	 
Let us consider the examples of the solution in simple particular cases. For two spatial variables and two equations~$K=2,n=2$
the matrices~$D_{x},D_{y},A $, satisfying all the conditions from~\cite{b1}, are given by 
\begin{equation}
	\label{matr1}
	D_{x}=\begin{pmatrix}
		d_{x1} & 0 \\ 0 & d_{x2}
	\end{pmatrix} ,
	D_{y}=\begin{pmatrix}
		d_{y1} & 0 \\ 0 & d_{y2}
	\end{pmatrix},
	A=\begin{pmatrix}
		-a & b \\ ka & -kb
	\end{pmatrix},
\end{equation}	
where $a>0,b>0,k>0$ are constants.

	Let us denote ${\displaystyle ab(a+bk)^{-2}}$ by $P,$ $d_{1x} - d_{2x}$ by $\Delta_x,$ $d_{1y} - d_{2y}$ by $\Delta_y.$ 
	
We compute 
	$M_{1} = -P\Delta_x^2$, $M_{2} =-P\Delta_x\Delta_y$, $M_{3}=-P\Delta_y^2$.
The matrix
\begin{equation}
	\label{matr2}
	\begin{pmatrix}
		M_{1} & M_{2} \\ M_{2} & M_{3}
	\end{pmatrix}
\end{equation}	
has the following eigenvalues
$\Lambda_{1}, \Lambda_{2}$:
\begin{align*}
	&\Lambda_{1} = 0;
	&\Lambda_{2} = (M_1+M_3)=-P (\Delta_x^2+\Delta_y^2) < 0.
\end{align*}

	For three spatial variables~$(K=3)$:	
	$V_1=(d_{1x},d_{1y},d_{1z})$,$V_2=(d_{2x},d_{2y},d_{2z})$, 
	$\bar\Delta=(\Delta_x,\Delta_y,\Delta_z)=\bar V_1-\bar V_2$, $ V_1 \neq V_2$.	
	 After some computations we obtain that in the case when~$K=3$  the matrix~$M$
	 has the eigenvalue~$0$ with the multiplicity~$2$ and one negative eigenvalue: 
	$\Lambda_1=-P(\Delta_x^2+\Delta_y^2+\Delta_z^2),\Lambda_2=0,\Lambda_3=0$. 
	
	The similar results were obtained in the cases~$K=4,5$. For $K=4$ the eigenvalues of the matrix~$M$
	are equal to~$\Lambda_1=-P(\Delta_1^2+\Delta_2^2+\Delta_3^2+\Delta_4^2),\Lambda_2=0,\Lambda_3=0,\Lambda_4=0,$ for $K=5:$   $\Lambda_1=-P(\Delta_1^2+\Delta_2^2+\Delta_3^2+\Delta_4^2+\Delta_5^2),\Lambda_2=\Lambda_3=\Lambda_4=\Lambda_5=0$. 
	
	For the greater number of equations~$n\geq3$ in the article~\cite{b1} numeric com\-pu\-ta\-ti\-ons for the number matrices~$D_i,A$ were performed. These computations allowed to formulate the hypothesis on the rank of the matrix~$M.$
	
	In the present work computations in the Wolfram Mathematica 11.3 computer algebra system lead us to the conjecture for Theorem~\ref{th:main}. It has to be noted that the computational complexity of finding of the eigenvalues increases drastically with a growth of number of equations and spatial variables. For example, in the case of three spatial variables we have to compute the eigenvalues for the $3\times3$ matrix (it is not presented here due to its large size), which contains the polynomial elements with~$15$ variables, consisting of~$400$ to~$450$ monomials. The direct approach to the solution in this case is impossible and to solve it we used the algorithm of finding and replacing the repeating expressions in elements of the matrix~$M.$ 
	
	 The main result of Theorem~\ref{th:main} allows us to do important conclusions on the spatial structure of the solution to the equation~(\ref{c11}), and the asymptotic expansion of the solution~(\ref{c6a}) to the problem~(\ref{c3}).
 

\end{document}